

\input amstex

\documentstyle{amsppt}

\loadbold

\magnification=\magstep1

\NoBlackBoxes

\pagewidth{6.5truein}
\pageheight{9.0truein}

\font\smallsl=cmsl8

\def\ad{\operatorname{ad}}

\def\max{\operatorname{max}}

\def\spec{\operatorname{spec}}

\def\prim{\operatorname{prim}}

\def\A{\Cal{A}}

\def\Fe{F_\epsilon}

\def\Fo{F_0}

\def\Fep{F_\epsilon^\varphi }

\def\Fop{F_0}

\def\refA{1}

\def\refAST{2}

\def\refC{3}

\def\refCVone{4}

\def\refCVtwo{5}

\def\refDKP{6}

\def\refDL{7}

\def\refDP{8}

\def\refGW{9}

\def\refJ{10}

\def\refLone{11}

\def\refLtwo{12}

\def\refLu{13}

\def\refMR{14}

\def\refMo{15}

\def\refMoS{16}

\def\refPW{17}

\def\refR{18}

\def\refRo{19}

\def\refW{20}


\topmatter



\pretitle{

\hbox{}

\vskip -.9truein

\hfill {\smallsl April 1998. To appear in Bulletin of the London
Mathematical Society.}

\vskip .25truein

         }


\title Noetherian Centralizing Hopf Algebra Extensions \\
and Finite Morphisms of Quantum Groups \endtitle

\rightheadtext{Noetherian Hopf Algebras}

\author Edward S. Letzter \endauthor

\address Department of Mathematics, Texas A{\&}M University, College
Station, TX 77843\endaddress

\email letzter\@math.tamu.edu\endemail

\abstract We study finite centralizing extensions $A \subset H$ of
Noetherian Hopf algebras. Our main results provide necessary and
sufficient conditions for the fibers of the surjection $\spec H
\twoheadrightarrow \spec A$ to coincide with the $X$-orbits in $\spec
H$, where $X$ denotes the finite group of characters of $H$ that
restrict to the counit of $A$. In particular, all of the fibers are
$X$-orbits if and only if the fiber over the augmentation ideal of $A$
is an $X$-orbit. An application to the representation theory of
quantum function algebras, at roots of unity, is
presented. \endabstract

\thanks This research was supported in part by grants from the
National Science Foundation. \endthanks

\subjclass 16D30, 16S20, 16P40, 16W30, 81R50 \endsubjclass

\endtopmatter

\document

\head 1. Introduction \endhead

If $A \subset H$ is an extension of finitely generated commutative
Hopf algebras (over an algebraically closed field $k$), then classical
duality provides an exact sequence of affine algebraic groups,
$$1 \rightarrow X \rightarrow \max H \rightarrow \max A \rightarrow
1,$$
where the multiplication of maximal ideals is obtained from the
convolution of their associated linear characters (i.e., algebra
homomorphisms onto $k$), and where $X$ consists of those linear
characters of $H$ that restrict to the counit on $A$. However,
supposing now that $A \subset H$ is a noncommutative extension of Hopf
algebras, the preceding group structure does not -- in most cases --
apply to their primitive spectra. But $X$ can still be defined as
before, and there are actions of $X$ on $H$ by automorphisms that fix
$A$ pointwise (as explained in (2.1) below). Moreover, these
automorphic actions on $H$ induce actions on its primitive spectrum,
and these actions on the primitive spectrum reduce to the left and
right multiplication by $X$ on $\max H$, under convolution, that
appeared in the original commutative case (see (3.1)). One can
therefore ask, in the noncommutative setting, whether the fibers of
the correspondence from $\spec H$ to $\spec A$ are determined by the
action of $X$, in a manner generalizing the classical theory. Our aim
in this paper is to present an affirmative answer to this question,
under hypotheses applicable to quantum function algebras at roots of
unity.

Now assume that $G$ is a connected and simply connected complex
semisimple Lie group, that $\ell$ is an odd positive integer prime to
$3$ if $G$ possesses a $G_2$ component, and that $\epsilon$ is a
primitive $\ell$th root of unity. Our main motivating example involves
the Hopf algebra embedding $\Fo \subset \Fe$, as studied by De
Concini, Lyubashenko, and Procesi \cite{\refDL, \refDP}
(multiparameter analogs may be found in the work of Costantini and
Varagnolo \cite{\refCVone, \refCVtwo}). Here, $\Fo$ is isomorphic to,
and $\Fe$ is a quantization of, the classical coordinate ring of
$G$. Moreover, $\Fo$ is contained within the center of $\Fe$, and
$\Fe$ is a finitely generated projective $\Fo$-module
\cite{\refDL}. (Following the terminology of \cite{\refPW, 1.8}, the
embedding of $\Fo$ into $\Fe$ may be viewed as a covering of quantum
groups.) Consequently, the group $X$ of characters of $\Fe$ that
restrict to the counit on $\Fe$ is finite, and it follows from
\cite{\refDP} that the fibers of the surjection $\prim \Fe
\twoheadrightarrow \prim \Fo$ are precisely the $X$-orbits in $\prim
\Fe$, under certain addtional restrictions on $\ell$. As an
application of our analysis below, the additional assumptions on
$\ell$ can be removed; see (3.2i). (This last result also follows from
recent independent work by Montgomery and Schneider \cite{\refMoS}
concerning Hopf Galois extensions.)

To briefly describe the primary context for our study, and our main
theorem, let $H$ be a noetherian Hopf algebra containing an
associative subalgebra $A$. Suppose further that $A$ is a right
coideal of $H$ and that $H$ is a finite centralizing extension of $A$
(i.e., $H$ is generated as an $A$-module by finitely many elements $x$
such that $ax = xa$ for all $a \in A$). As before, let $X$ denote the
(finite) set of characters of $H$ that restrict to the counit on
$A$. In (2.9) we conclude that the following conditions are
equivalent: (a) the primitive ideals of $H$ contracting to the
augmentation ideal of $A$ all have codimension $1$ in $H$, (b) the
primitive ideals of $H$ contracting to the augmentation ideal of $A$
comprise a single $X$-orbit in $\prim H$, (c) the fibers of the
surjection $\prim H \twoheadrightarrow \prim A$ are precisely the
$X$-orbits in $\prim H$, (d) the fibers of the surjection $\spec H
\twoheadrightarrow \spec A$ are precisely the $X$-orbits in $\spec H$.

The material in this paper was presented in the special session on
Algebras, Algebra Cohomology, and Polynomial Identities, at the 102nd
Annual Meeting of the American Mathematical Society, Orlando, January,
1996.

I am grateful to M. Costantini for the helpful comments used in (3.2).

\head 2. Noetherian Centralizing Hopf Algebra Extensions \endhead

The main results of this section are found in (2.7) and (2.9).

\subhead 2.1 Preliminaries \endsubhead All of the algebras considered
in this section are defined over a single base field $k$.

(i) The comultiplication and counit (augmentation map) of a bialgebra
will always be denoted $\Delta$ and $\varepsilon$, respectively, and the
antipode of a Hopf algebra will always be denoted $S$. We will use the
notation $\Delta (b) = \sum b_1 \otimes b_2$.

(ii) Let $B$ be a bialgebra, and let $\chi$ be a (linear) character of
$B$ (i.e., $\chi$ is an algebra homomorphism of $B$ onto $k$). Let
$k_\chi$ denote the one-dimensional $B$-module defined by the action
$b.1 = \chi(b)$. The assignment
$$\sigma _\chi\colon b \mapsto \sum \chi(b_1)b_2, \qquad \text{for $b
\in B$,}$$
is an algebra endomorphism, which we will call, following
\cite{\refJ}, a {\sl (right) winding\/} endomorphism of $B$. If $B$ is
a Hopf algebra, then every winding endomorphism is an automorphism;
see, e.g., \cite{\refJ, 1.3.4}. Furthermore, if $\chi'$ is also a
character of $B$, then $\sigma _\chi \circ \sigma _{\chi'} = \sigma
_{\chi' \ast \chi}$, where $\chi' \ast \chi$ is the convolution of
these characters. Consequently, if $B$ is a Hopf algebra then the
character group of $B$ acts on $B$ by automorphisms.

(iii) Retaining the preceding notation, let $A$ be an associative
subalgebra of $B$. Observe, when $\sigma_\chi$ acts as the identity on
$A$, that $\chi \left\vert _A \right. = \varepsilon \left\vert _A
\right.$, because
$$\chi(a) \; = \; \chi\left(\sum a_1\varepsilon(a_2)\right) \; = \;
\sum \chi (a_1)\varepsilon (a_2) \; = \; \varepsilon \left(\sum
\chi(a_1)a_2 \right),$$
for all $a \in A$. Now assume that $A$ is a right coideal of $B$
(i.e., $\Delta(A) \subseteq A\otimes B$). In converse to the
preceding, if $\chi$ acts as $\varepsilon$ on $A$ then $\sigma_\chi$
restricts to the identity on $A$. Also, the set $X$ of characters of
$B$ restricting to $\varepsilon$ on $A$ is closed under convolution,
and so forms of subgroup of the character group of $B$.

\subhead 2.2 \endsubhead Let $H$ be a Hopf algebra, and let $V$ be an
$H$-$H$-bimodule. Recall that the $\ad H$-module structure on $V$ is
defined by the left action
$$\ad (h) v \; = \; \sum h_1.v.S(h_2),$$
where $h \in H$, $v \in V$, and $\Delta(h) = \sum h_1 \otimes
h_2$. 

Further recall, when $V$ has a one-dimensional $\ad H$-submodule
$k$-spanned by the element $n$, that
$$\multline n.\sigma _\chi(h) \; = \; \sum n.\chi (h_1)h_2 \; = \;
\sum (\ad (h_1)n).h_2 \; = \\ \sum h_1.n.S(h_2)h_3 \; = \; \sum
h_1.n.\varepsilon(h_2) \; = \; \sum h_1\varepsilon(h_2).n \; = \; h.n,
\endmultline$$
for all $h \in H$, where $\chi$ is the character of $H$ defined by
$\ad (h)n = \chi(h)n$.

\subhead 2.3 \endsubhead We now review some elementary facts from
noetherian ring theory (the reader is referred to \cite{\refGW} or
\cite{\refMR} for further details). Let $R_\alpha$ and $R_\beta$ be
prime noetherian rings, and let $M$ be an
$R_\alpha$-$R_\beta$-bimodule that is finitely generated on each side.

(i) Suppose that $M$ is torsionfree (cf\. \cite{\refMR, 3.4.2}) on
each side; $M$ is then said to be a {\sl bond\/} from $R_\alpha$ to
$R_\beta$. Next, let $m$ be a nonzero normal element of $M$ (i.e.,
$R_\alpha.m = m.R_\beta$), and let $a$ be an element of $R_\alpha$ for
which $a.m = 0$. Then $0 = a.m = R_\alpha a.m.R_\beta = R_\alpha
aR_\alpha.m$, and so $a = 0$. (Recall, from Goldie's Theorem, that
every nonzero ideal in a noetherian prime ring contains a regular
element.) We see, therefore, that the left and right annihilators of
any nonzero normal element of $M$ are both equal to zero. If
$R_\gamma$ is also a prime noetherian ring and $N$ is a bond from
$R_\beta$ to $R_\gamma$, then there exists an
$R_\alpha$-$R_\gamma$-bimodule factor of $M\otimes _{R_\beta}N$ that
is a bond from $R_\alpha$ to $R_\gamma$; see, e.g., \cite{\refW, 5.1}.

(ii) Let $R$ be a noetherian ring, and let $N$ be an $R$-$R$-bimodule
that is finitely generated on each side. Let $P_\alpha$ and $P_\beta$
be prime ideals of $R$, and suppose that $N$ is a bond from
$R/P_\alpha$ to $R/P_\beta$. Assume further that there exists an
automorphism $\sigma$ of $R$, and a nonzero element $n$ of $N$, such
that $r.n = n.\sigma(n)$ for all $r \in R$. Observe that
$0 = P_\alpha.n = n.\sigma(P_\alpha)$ and that $0 = n.P_\beta =
\sigma^{-1}(P_\beta).n$. Consequently, using (i),
$\sigma(P_\alpha) = P_\beta$.

\subhead 2.4 \endsubhead Let $H$ be a noetherian Hopf algebra
containing a noetherian subalgebra $A$ such that $H$ is finitely
generated as a left and right $A$-module. Assume further that every
$H$-module composition factor of the finite dimensional left
$H$-module $H\otimes_A k_\varepsilon$ is one dimensional, and let $X$
denote the set of characters of $H$ corresponding to these composition
factors.

\proclaim{Lemma} Let $P_\alpha$ and $P_\beta$ be prime ideals of $H$,
and let $Q$ be a prime ideal of $A$.

{\rm (i)} Suppose that there exists an $H$-$A$-bimodule factor of $H$
bonding $H/P_\alpha$ to $A/Q$ and an $A$-$H$-bimodule factor of $H$
bonding $A/Q$ to $H/P_\beta$. Then $\sigma_{\chi}(P_\alpha) = P_\beta$
for some $\chi \in X$.

{\rm (ii)} Suppose that $P_\alpha \cap A$ and $P_\beta \cap A$ are
both equal to $Q$. Then $\sigma _{\chi}(P_\alpha) = P_\beta$ for some
$\chi \in X$.  \endproclaim

\demo{Proof} (i) By (2.3i), there is a nonzero $H$-$H$-bimodule factor
$V$ of $H \otimes _A H$ that is a bond from $H/P_\alpha$ to
$H/P_\beta$. Let $e$ denote the image in $V$ of $1\otimes 1$. Observe
that $e$ is nonzero since it is the generator of $V$ as an
$H$-$H$-bimodule. Furthermore, since $a.e = e.a$ for all $a \in A$, it
follows that the $\ad H$-module generated by $e$ is isomorphic, as a
left $H$-module, to a nonzero factor of $H \otimes _A
k_\varepsilon$. Consequently, $V$ contains an element $n$ such that $\ad
(h)n = \chi(h)n$ for some $\chi
\in X$. It now follows from (2.2) that $h.n = n.\sigma_\chi(h)$, for
all $h \in H$, and it follows from (2.3ii) that $\sigma_\chi(P_\alpha) =
P_\beta$.

(ii) It follows, for example, from \cite{\refLone, 1.1} that there
exists an $H$-$A$-bimodule factor of $H/P_\alpha$ bonding $H/P_\alpha$
to $A/Q$. Similarly, there exists an $A$-$H$-bimodule factor of
$H/P_\beta$ bonding $A/Q$ to $H/P_\beta$. Part (ii) now follows from
part (i). \qed\enddemo

\subhead 2.5 \endsubhead (i) Let $R$ be an algebra. The set of prime
ideals of $R$ will be denoted $\spec R$ and the set of (left)
primitive ideals of $R$ will be denoted $\prim R$; these sets will be
endowed with the Jacobson (Zariski) topology. If $\A$ is a group
acting on $R$ by automorphisms, one immediately obtains actions of
$\A$ on $\spec R$ and $\prim R$.

(ii) Assume that $R$ is a {\sl finite centralizing\/} extension of a
subalgebra $U$ (i.e., assume that $R = Ur_0 + \cdots + Ur_s$, for
elements $r_0,\ldots,r_s \in R$, such that $ur_i = r_iu$ for all $0
\leq i \leq s$ and all $u \in U$). It is well known that $U$ is
noetherian if and only if $R$ is noetherian (see, e.g., \cite{\refMR,
10.1.11}) and that if $P$ is a prime ideal of $U$ then $P \cap U$ is a
prime ideal of $R$ (e.g., \cite{\refMR, 10.2.4}).

(iii) Assigning $P\cap U$ to each prime ideal $P$ of $R$ produces
closed continuous surjections
$$u\colon\spec R \twoheadrightarrow \spec U \quad \text{and} \quad
u\colon\prim R \twoheadrightarrow \prim U ,$$
with finite fibers. See, for example, \cite{\refMR, Chapter 10} for
more details. (Finiteness of the fibers follows from \cite{\refRo,
3.4}.)

(iv) Let $M$ be a simple left $U$-module, and let $r_0,\ldots,r_s$ be
as in (ii). Then $R\otimes _UM = r_0\otimes M + \cdots + 
r_s\otimes M$, and we see that $R\otimes _UM$ is a semisimple
$U$-module, of length no greater than $s$, with each $U$-module
composition factor isomorphic to $M$.

(v) We will say that an $R$-module $V$ {\sl restricts\/} to a
character $\xi$ of $U$ if $u.v = \xi (u)v$ for all $u \in U$ and $v
\in V$. Observe that any finitely generated $R$-module restricting to
a character of $U$ must be finite dimensional. Furthermore, it follows
from (iv) that a simple $R$-module $W$ restricts to $\xi$ if and only
if $W$ is an $R$-module composition factor of $R\otimes
_Uk_\xi$. Consequently, there are only finitely many simple
$R$-modules restricting to $\xi$, and so there are at most finitely
many characters of $R$ restricting to $\xi$. (This last assertion also
follows from (iii).)

\subhead 2.6 \endsubhead Let $B$ be a bialgebra that is a finite
centralizing extension of an associative subalgebra $A$. 

(i) In view of (2.5v), the set of characters of $B$ restricting to
$\varepsilon$ on $A$ is finite.

(ii) The kernel of $\varepsilon \left \vert _A \right.$ will be
denoted $A^+$. Note that $BA^+ = A^+B$, that $BA^+$ is contained
within the augmentation ideal of $B$, and that $B/BA^+$ is a finite
dimensional algebra. Of course, if $A$ is a Hopf algebra, and if the
inclusion of $A$ in $B$ is a homomorphism of Hopf algebras, then
$B/BA^+$ is the Hopf cokernel.

\subhead 2.7 \endsubhead We now record our primary abstract result,
stated in its most general form. A more condensed version is presented
in (2.9).

\proclaim{Theorem} Let $H$ be a noetherian Hopf algebra that is
a finite centralizing extension of an associative subalgebra $A$, and
let $X$ denote the set of characters of $H$ that restrict to
$\varepsilon$ on $A$. The following conditions are equivalent:

\noindent {\rm (i)} Every irreducible $H/HA^+$-module is one
dimensional.

\noindent {\rm (ii)} Every irreducible $H$-module restricting to
$\varepsilon$ on $A$ is one dimensional.

\noindent {\rm (iii)} If $P_\alpha$ and $P_\beta$ are primitive ideals
of $H$ for which $P_\alpha \cap A = P_\beta \cap A$, then there exists
$\chi \in X$ such that $\sigma_\chi (P_\alpha) = P_\beta$.

\noindent {\rm (iv)} If $P_\alpha$ and $P_\beta$ are prime ideals
of $H$ for which $P_\alpha \cap A = P_\beta \cap A$, then there exists
$\chi \in X$ such that $\sigma_\chi (P_\alpha) = P_\beta$.
\endproclaim

\demo{Proof} $\roman{(ii)} \Rightarrow \roman{(iv)}$ Let $P_\alpha$
and $P_\beta$ be prime ideals of $H$ for which $P_\alpha\cap A =
P_\beta \cap A$.  It follows from (2.5v) and (ii) that every
$H$-module composition factor of $H \otimes _A k_\varepsilon$ is one
dimensional. Statement (iv) now follows from (2.4ii).

The remaining implications are immediate. \qed\enddemo

\subhead 2.8 \endsubhead Retaining the notation of (2.7), assume
further that $A$ is a right coideal of $H$. It follows from (2.1iii)
and (2.5v) that $X$ is a finite subgroup of the character group of
$H$, and that the set $\left\{ \sigma_\chi \mid \chi \in X\right\}$ is
precisely the group of winding automorphisms of $H$ that fix $A$
pointwise. Consequently, each fiber of the surjection $u\colon \spec H
\twoheadrightarrow \spec A$ is a union of $X$-orbits.

The following is now an immediate corollary to (2.7), and presents a
somewhat smoother alternative to that result.

\proclaim{2.9 Theorem} Let $H$ be a noetherian Hopf algebra that is a
finite centralizing extension of an associative subalgebra $A$, and
suppose that $A$ is a right coideal of $H$. Let $X$ be the finite
group of characters of $H$ that restrict to the counit on $A$.

The following conditions are equivalent:

{\rm (i)} Every irreducible $H$-module restricting to
the counit on $A$ is one dimensional.

{\rm (ii)} The primitive ideals of $H$ contracting to the augmentation
ideal of $A$ comprise a single $X$-orbit in $\prim H$.

{\rm (iii)} The fibers of the surjection $u\colon \prim H \twoheadrightarrow
\prim A$ are precisely the $X$-orbits in $\prim H$.

{\rm (iv)} The fibers of the surjection $u\colon \spec H \twoheadrightarrow
\spec A$ are precisely the $X$-orbits in $\spec H$.  \endproclaim

\remark{2.10 Remark} Retaining the assumptions of (2.9), but
additionally supposing that $A$ is a Hopf subalgebra of $H$, it turns
out that the choice of $\varepsilon$ in (2.9i) is somewhat arbitrary. To
explain, let $\xi$ be any character of $A$ for which every irreducible
$H$-module restricting to $\xi$ is one dimensional, and let $\chi$
denote one of the characters of $H$ restricting to $\xi$. Set $\sigma
= \sigma_\chi$. Because $A$ is a Hopf subalgebra, $\sigma$ maps $A$ to
itself automorphically. Now let $K$ be the kernel of $\xi$, and
observe for $a \in K$ that $\varepsilon(\sigma (a)) = \sum \chi
(a_1)\varepsilon(a_2) = \chi (a) = \xi (a) = 0$. Therefore, $\sigma (K) =
A^+$, and so $\sigma (HK) = HA^+$. In particular, $\sigma$ induces an
algebra isomorphism from $H/HK$ onto $H/HA^+$, and so the irreducible
$H$-modules restricting to $\varepsilon$ are also all one dimensional.
\endremark

\head 3. Finite Morphisms of Quantum Groups \endhead

In this part we outline some connections between the results in the
previous section and the representation theory of quantum groups.

\subhead 3.1 \endsubhead We first review the classical case
(cf\. \cite{\refA, Chapter 4; \refMo, \S 9.3}), which was briefly
sketched in the introduction. To begin, assume that $k$ is an
algebraically closed field, and let
$$1 \rightarrow X \rightarrow Y @>\pi>> Z \rightarrow 1$$
be a short exact sequence of $k$-affine algebraic groups. Recall that
there is a corresponding embedding $A \subset H$ of finitely generated
reduced commutative Hopf algebras, where $Z$ is the character group
of $A$, and where $Y$ is the character group of $H$. Conversely, an
epimorphism of algebraic groups can be associated to each embedding of
finitely generated commutative $k$-Hopf algebras (that need not be
assumed reduced -- see, e.g., \cite{\refMo, 9.2.11--12}).

Observe that $X$, viewed as a subgroup of $Y$, is equal to the group
of characters of $H$ that act as $\varepsilon$ on $A$. Hence, by
(2.1iii), $X$ acts on $H$ by winding automorphisms, and $\{ \sigma
_\gamma \mid \gamma \in X \}$ is precisely the set of winding
automorphisms of $H$ that fix $A$ pointwise.

On the other hand, by imposing the group structures of $Y$ and $Z$
(i.e., convolution of characters) onto $\max H$ and $\max A$,
respectively, we see that the subgroup $X$ of $Y$ acts on $\max H$
by left and right multiplication. Moreover, since $X$ is the kernel
of $\pi$, it immediately follows that the induced multiplication
actions on $\max A$ are trivial and that the fibers of the surjection
$\max H \twoheadrightarrow \max A$ are exactly the $X$-orbits, under
the right or left multiplication action, in $\max H$. Finally, it is
easy to check that the $X$-orbits in $\max H$ under the
multiplication actions are precisely the $X$-orbits, in $\max H$,
induced by the winding automorphisms mentioned in the preceding
paragraph. 

In particular, for finite extensions of finitely generated commutative
Hopf algebras over algebraically closed fields, the equivalent
conditions in (2.9) may be deduced from the classical theory.

\subhead 3.2 \endsubhead Now let $G$ be a connected, simply connected,
complex semisimple Lie group. Fix an odd positive integer $\ell$
(prime to $3$ if $G$ has a $G_2$ component), and let $\epsilon$ be a
primitive $\ell$th root of $1$.

(i) It is shown by De Concini and Lyubashenko that the complex quantum
function algebra $\Fe$ of $G$ (as specified in \cite{\refDL, \S 9})
contains a central sub-Hopf-algebra $\Fo$ isomorphic to the ring of
regular functions on $G$, and that $\Fe$ is finitely generated as an
$\Fo$-module \cite{\refDL, 6.4, 7.2}. Furthermore, it follows from
\cite{\refDL, 10.7} that every irreducible $\Fe$-module restricting to
the counit on $\Fo$ is one dimensional. Hence, by (2.9), the fibers of
the surjections
$$\spec \Fe \twoheadrightarrow \spec \Fo \quad \text{and} \quad \prim
\Fe \twoheadrightarrow \prim \Fo $$
are precisely the $X$-orbits in $\spec \Fe$ and $\prim \Fe$,
respectively, where $X$ is the group of characters of $\Fe$
restricting to the counit of $\Fo$.

The preceding conclusion has been independently verified by
S. Montgomery and H.-J. Schneider, in their study of Hopf Galois
extensions \cite{\refMoS}. Also, under the additional assumption that
$\ell$ is prime to the bad primes of the associated root system, the
conclusion already follows from \cite{\refDP, \S 4.10}. The
approach in \cite{\refDL; \refDP} involves a detailed analysis of
certain Azumaya algebras, obtained as localizations of factors of
$\Fe$, and precise calculations of the dimensions of the irreducible
representations are obtained.

(ii) Assuming that $G$ is simple, an extension $\Fop \subset \Fep$ of
complex Hopf algebras is studied by Costantini and Varagnolo
\cite{\refCVone; \refCVtwo}, where $\Fep$ is a multiparametric
quantization of the coordinate algebra of $G$ (cf\. \cite{\refDKP;
\refR}) and $\Fop$ is isomorphic to the ring of regular functions on
$G$. It is shown in \cite{\refCVone} that $\Fop$ is contained within
the center of $\Fep$ and that $\Fep$ is finitely generated as an
$\Fop$-module.  Furthermore, it follows as in \cite{\refDL} that every
irreducible representation of $\Fep$ restricting to the counit of
$\Fop$ is one dimensional \cite{\refC}, and so the conditions in (2.9)
hold for $\spec \Fep \twoheadrightarrow \spec\Fop$. Under certain
additional restrictions on $\ell$, this result follows from
\cite{\refCVtwo}, where dimensions of the irreducible representations
are also calculated.

(iii) In contrast, it is a fundamental property of the central Hopf
algebra embeddings featured in the representation theory, at roots of
unity, of quantized enveloping algebras of semisimple Lie algebras
(see, e.g., \cite{\refDP, Chapter 5; \refLu, Chapter 35}) that the
conditions in (2.9) are not satisfied.

(iv) Letting $k$ be an arbitrary field, there is a $k$-bialgebra
embedding of the classical coordinate ring $k[M_n]$ of $n\times n$
matrices into the center of the quantum coordinate ring $k_q[M_n]$ of
$n\times n$ matrices, where $q$ is a primitive $t$th root of unity in
$k$ and $t$ is odd (cf\. \cite{\refPW, Chapter 7}). Letting $X$ denote
the group of convolution invertible characters of $k_q[M_n]$ that
restrict to the counit on $k[M_n]$, there is an action of $X \times X$
on $k_q[M_n]$, by right and left winding automorphisms, that fixes
$k[M_n]$ pointwise. It is shown in \cite{\refLtwo, 2.12} that the fibers
of $\spec k_q[M_n] \twoheadrightarrow \spec k[M_n]$ coincide with the
$X\times X$ orbits in $\spec k_q[M_n]$. 

 \enddocument